\documentclass{article}

\usepackage[utf8]{inputenc}
\usepackage[english]{babel}
\usepackage{geometry}
\usepackage{hyperref}
\hypersetup{
	colorlinks=true,
	linkcolor=blue,
	urlcolor=blue,
	filecolor=blue,
	citecolor=red,
}
\usepackage{amssymb}
\usepackage{amsthm}
\usepackage{amsmath}
\usepackage{mathrsfs}
\usepackage{esint}
\usepackage{stix}
\usepackage{authblk}
\usepackage{blindtext}

\theoremstyle{plain}
\newtheorem{theorem}{Theorem}[section]
\newtheorem{lemma}[theorem]{Lemma}
\newtheorem{corollary}[theorem]{Corollary}

\theoremstyle{definition}
\newtheorem{definition}[theorem]{Definition}

\theoremstyle{remark}
\newtheorem{remark}[theorem]{Remark}

\numberwithin{equation}{section}

\begin{document}
	
	\title{A note on Bernoulli type free boundary problem on collapsed $RCD(K,N)$-spaces}
	
	\author[1]{Sitan Lin}
	\affil[1]{School of Mathematical Sciences, and Institute of Natural Sciences,\authorcr Shanghai Jiao Tong University,\authorcr e-mail: lin\_sitan@sjtu.edu.cn}
	
	\date{\today}
	
	\maketitle
	
\begin{abstract}
	In this paper, we investigate Bernoulli type free boundary problem on collapsed $RCD(K,N)$-spaces. We prove the existence of minimizers and prove the local Lipschitz continuity of minimizers provided that the negative part is locally Lipschitz continuous. In particular, we prove the local Lipschitz continuity of minimizers for the one-phase problem (i.e. when the solution is non-negative). And then we prove that the free boundaries of minimizers have locally finite perimeter. We emphasize that the proof in this paper applies to collapsed $RCD(K,N)$-spaces and does not rely on the non-collapsed condition.
\end{abstract}
	
\section{Introduction}
	The Bernoulli problem is among the most extensively studied free boundary problems on $\mathbb{R}
	^{n}$. It generally serves as one of the most fundamental free boundary models. Beyond its physical motivation from fluid dynamics, the Bernoulli problem has mathematical connections to minimal surfaces and shape optimization problems. An elementary formulation of the (two-phase) Bernoulli problem is as below: a continuous function $u$ on a domain $\Omega$ such that
	\begin{equation}\label{equation_1_Introduction}
		\begin{cases}
			\Delta u_{\pm} =0\quad & \text{in } \{u_{\pm}>0\}\cap\Omega,\\
			u_{\pm}=0\quad & \text{on } \partial\{u_{\pm}>0\}\cap\Omega,\\
			\partial_{\nu}u_{\pm}=\sqrt{\lambda_{\pm}}\quad & \text{on } \left(\partial\{u_{\pm}>0\}\setminus\partial\{u_{\mp}>0\}\right)\cap\Omega,\\
			(\partial_{\nu}u_{+})^{2}-(\partial_{\nu}u_{-})^{2}=\lambda_{+}-\lambda_{-}\quad & \text{on } \left(\partial\{u_{+}>0\}\cap\partial\{u_{-}>0\}\right)\cap\Omega,
		\end{cases}
	\end{equation}
	where $u_{+}$ and $u_{-}$ are the positive part and negative part of $u$ respectively, $\nu$ is the inward normal to the set $\{u_{\pm}>0\}$ and $\lambda_{+}>\lambda_{-}>0$ are constants. Here, $\{u_{\pm}>0\}$ is presumed to be sufficiently smooth. The problem is overdetermined if the sets $\{u_{\pm}>0\}$ are given, but since the sets $\{u_{\pm}>0\}$ are free, it is possible to find a solution. The problem \eqref{equation_1_Introduction} is said to be one-phase if $u$ is non-negative.
	
	Let us consider the energy functional
	\begin{equation}\label{equation_2_Introduction}
		E(u):=\int_{\Omega} \lvert\nabla u\rvert^{2} + \lambda_{+}\chi_{\{u^{+}>0\}} + \lambda_{-}\chi_{\{u^{-}>0\}}.
	\end{equation}
	If $u$ is a critical point of $E$, then it is easily checked that $u$ is harmonic in $\{u_{\pm}>0\}$. If, in addition, $\Omega$ is a domain in $\mathbb{R}^{n}$ and the free boundaries $\partial\{u_{\pm}>0\}\cap\Omega$ are sufficiently smooth, then the domain variation gives that
	\begin{equation}
		\int_{\partial\{u_{+}>0\}} \left(\lvert\nabla u_{+}\rvert^{2}-\lambda_{+}\right)\xi\cdot\nu\mathrm{d}\mathcal{H}^{n-1} + \int_{\partial\{u_{-}>0\}} \left(\lvert\nabla u_{-}\rvert^{2}-\lambda_{-}\right)\xi\cdot\nu\mathrm{d}\mathcal{H}^{n-1}=0
	\end{equation}
	for each smooth vector field $\xi\in C_{0}^{\infty}(\Omega;\mathbb{R}^{n})$. In particular, $u$ is a classical solution to \eqref{equation_1_Introduction}.
	
	Since minimizers of \eqref{equation_2_Introduction} are critical points, a natural way to find a solution of the Bernoulli problem \eqref{equation_1_Introduction} is to minimize the energy functional \eqref{equation_2_Introduction}. Following \cite{Alt-Caffarelli_1981}, a regularity theory for minimizers on $\mathbb{R}^{n}$ is by now well-developed. Roughly speaking, minimizers of \eqref{equation_2_Introduction} are locally Lipschitz continuous \cite{Alt-Caffarelli_1981,Alt-Caffarelli-Friedman_1984} and the free boundaries are smooth hypersurfaces \cite{Alt-Caffarelli_1981,De_Philippis-Spolaor-Velichkov_2021,Kinderlehrer-Nirenberg_1977}, except on a set of Hausdorff dimension $\le n-k^{*}$ with $k^{*}\in\{5,6,7\}$ \cite{Weiss_1999,Caffarelli-Jerison-Kenig_2004,De_Silva-Jerison_2009,Jerison-Savin_2015}.
	
	In this paper, we investigate the Bernoulli problem on $RCD(K,N)$-spaces. Roughly speaking, $RCD(K,N)$-spaces are metric measure spaces with Ricci curvature bounded from below by $K\in\mathbb{R}$ and dimension bounded from above by $N\in[0,\infty]$ in the synthetic sense. Examples of $RCD(K,N)$-spaces includes all smooth Riemannian manifolds with Ricci curvature $\ge K$ and dimension $\le N$. It is worth mentioning that the class of $RCD(K,N)$-spaces is closed under measured Gromov-Hausdorff convergence, while that of smooth Riemannian manifolds is not. We refer the readers to the survey \cite{Ambrosio_2018} and the references therein for backgrounds on $RCD(K,N)$-spaces.
	
	In this paper, we focus on the absolute minimizers of the energy functional \eqref{equation_2_Introduction}. The reasons are as follows: (1) We lack the first variation formula with respect to the domain variation for the volume term $\int_{\Omega} \chi_{\{u>0\}}\mathrm{d}m$ on $RCD(K,N)$-spaces. Consequently, we are unable to define the critical points (with respect to the domain variation) of the energy functional \eqref{equation_2_Introduction} on $RCD(K,N)$-spaces; (2) Unlike the Euclidean spaces, $RCD(K,N)$-spaces generally exhibit poor regularity (collapsed $RCD(K,N)$-spaces may even lack a manifold structure \cite{Hupp-Naber-Wang_2025,Zhou_2024}), which prevents us from discussing the classical solutions of equation \eqref{equation_1_Introduction} on general $RCD(K,N)$-spaces at the current stage. On the other hand, the first order calculus on metric measure spaces has been well-developed \cite{Cheeger_1999} (in particular, the weak compactness for the $W^{1,2}$-spaces and the lower semi-continuity for the Dirichlet energy with respect to the weak convergence have been established), which allows us to investigate the absolute minimizers of the energy functional \eqref{equation_2_Introduction} on $RCD(K,N)$-spaces.
	
	In \cite{Chan-Zhang-Zhu_2021}, the authors studied the one-phase Bernoulli problem on non-collapsed $RCD(K,N)$-spaces. They proved the existence and the local Lipschitz continuity of a minimizer. Moreover, under the assumption that $\Omega$ is disjoint from the boundary of the ambient space, they proved that the free boundary is a topological manifold of Hausdorff dimension $N-1$, except on a relatively closed subset of codimension two. We emphasize that the proof in \cite{Chan-Zhang-Zhu_2021} requires that the ambient space is non-collapsing.

	The main result of this paper is as follows:
	\begin{theorem}\label{theorem_1_Introduction}
		Let $(X,d,m)$ be an $RCD(K,N)$-space with $K\in\mathbb{R}$ and $N\in(1,\infty)$ and let $\Omega\subset X$ be a bounded open domain. Given $g\in W^{1,2}(\Omega)$, we define $\mathscr{A}_{g}:=\{v\in W^{1,2}(\Omega):v-g\in W^{1,2}_{0}(\Omega)\}$ and consider the following minimization problem:
		\begin{equation}
			\min_{v\in\mathscr{A}_{g}} E(v),
		\end{equation}
		where $E(v)=\int_{\Omega}\lvert\nabla v\rvert^{2}+\lambda_{+}\chi_{\{v>0\}}+\lambda_{-}\chi_{\{v<0\}}\mathrm{d}m$. Then the following statements hold.
		\begin{itemize}
			\item[(1)] There exists a minimizer $u\in\mathscr{A}_{g}$, i.e.
			\begin{equation}
				E(u)\le E(v),\quad\forall v\in\mathscr{A}_{g}.
			\end{equation}
			\item[(2)] A minimizer $u$ is locally Lipschitz continuous if and only if one of $u_{+}$ and $u_{-}$ is Lipschitz continuous. In particular, for the one-phase problem (i.e. if $u\ge 0$), $u$ is locally Lipschitz continuous.
			\item[(3)] The free boundaries $\partial\{u_{\pm}>0\}\cap\Omega$ of a minimizer $u$ have locally finite perimeter in $\Omega$.
		\end{itemize}
	\end{theorem}
	\begin{remark}\label{remark_1_Introduction}
		We also include a proof of the non-degeneracy for minimizers.
	\end{remark}
	\begin{remark}\label{remark_2_Introduction}
		Regarding the statement (2) of Theorem \ref{theorem_1_Introduction}, we point out that the proof in \cite{Chan-Zhang-Zhu_2021} relies on the assumption that $(X,d,m)$ is non-collapsed, i.e. $m$ is the $N$-dimensional Hausdorff measure and $N$ is an integer. The proof in this paper originates from \cite{Danielli-Petrosyan_2005}, where the authors dealt with $p$-Laplace operator on $\mathbb{R}^{n}$. Compared to \cite{Danielli-Petrosyan_2005}, we need to consider Sobolev functions defined on varying metric measure spaces. Consequently, we require the stability results from \cite{Ambrosio-Honda_2017}. Moreover, in \cite{Danielli-Petrosyan_2005}, the authors used the property that minimizers converge to minimizers, which plays a key role in the classical theory. However, on collapsed $RCD(K,N)$-spaces, it is not even known whether this property holds. In this paper, we avoid using this property and instead directly prove that the blow-up limit (of the sequence $u_{k}$ constructed in the proof of Lemma \ref{lemma_3_Lipschitz_continuity}) is a harmonic function by using the Laplacian estimate (Lemma \ref{lemma_2_Lipschitz_continuity}).
	\end{remark}
	\begin{remark}\label{remark_3_Introduction}
		Regarding the statement (3) of Theorem \ref{theorem_1_Introduction}, we point out that for the one-phase problem, if $(X,d,m)$ is non-collapsed, then the free boundary is a topological manifold except on a set of codimension two \cite{Chan-Zhang-Zhu_2021}. The proof in \cite{Chan-Zhang-Zhu_2021} relies on Cheeger-Colding's Reifenberg theorem \cite[Appendix 1]{Cheeger_1997}, which does not hold on collapsed $RCD$-spaces. In fact, there exists a collapsed $RCD(K,N)$-space $(Y,d_{Y},m_{Y})$ that does not have a manifold structure \cite{Hupp-Naber-Wang_2025}, see also \cite{Zhou_2024}. More precisely, any open subset $U\subset Y$ is not a topological manifold. On the other hand, $Y\times\mathbb{R}$ is an $RCD(K,N+1)$-space. We define $u:Y\times\mathbb{R}\to\mathbb{R}$ by
		\begin{equation}
			u(y,t):=t_{+},\quad y\in Y,t\in\mathbb{R}.
		\end{equation}
		Then it can be easily checked that $u$ is a minimizer while its free boundary is $Y$. This example demonstrates that, on general $RCD(K,N)$-spaces (possibly collapsed), the free boundary of a minimizer may not have a manifold structure.
		
		Recall that the reduced boundary of a set of locally finite perimeter is rectifiable \cite{Brue-Pasqualetto-Semola_2023}. In particular, the statement (3) of Theorem \ref{theorem_1_Introduction} implies that the free boundaries of minimizers are rectifiable.
	\end{remark}
	\begin{remark}\label{remark_4_Introduction}
		Note that there exists a non-minimizing solution to the Bernoulli problem \eqref{equation_1_Introduction} on $\mathbb{R}^{n}$. Indeed, it can be checked that the function
		\begin{equation}
			u(x):=\left(-\ln\lvert x-e_{1}\rvert\right)_{+} + \left(-\ln\lvert x+e_{1}\rvert\right)_{+}
		\end{equation}
		is a classical solution on $\Omega:=\mathbb{R}^{2}\setminus \left(B_{1/2}(e_{1})\cup B_{1/2}(-e_{1})\right)$, but it is not a minimizer. Indeed, $u$ is a variational solution defined in \cite{Weiss_2003,Kriventsov-Weiss_2025}. Roughly speaking, variational solutions are critical points of the energy functional \eqref{equation_2_Introduction} under domain variations. On general $RCD(K,N)$-spaces, although we can use the regular Lagrangian flow to derive the first variation formula with respect to domain variations for the Dirichlet term $\int_{\Omega} \lvert\nabla u\rvert^{2}\mathrm{d}m$, difficulties remain in establishing the first variation formula with respect to domain variations for the volume term $\int_{\Omega} \chi_{\{u>0\}}\mathrm{d}m$. Therefore, we focus on minimizers in this paper. The investigation of non-minimizing solutions is left for future research.
	\end{remark}
	
\section{Preliminaries}
In this section, we review some basic knowledge about $RCD$-spaces that will be used later in this paper. For more backgrounds on $RCD$-spaces, we refer the readers to the survey \cite{Ambrosio_2018} and the references therein.

Throughout this paper, we always denote by $(X,d,m)$ a metric measure space, i.e. $(X,d)$ is a complete and separable metric space and $m$ is a locally finite Borel measure on $X$. We also assume that $(X,d,m)$ is an $RCD(K,N)$-space with $K\in\mathbb{R}$ and $N\in(1,\infty)$.

\subsection{Calculus tools}
We denote by $\mathop{\mathrm{Lip}}(X)$ the collection of all Lipschitz functions on $(X,d)$. For each $f\in\mathop{\mathrm{Lip}}(X)$, the pointwise Lipschitz constant of $f$ at $x$ is defined by
\begin{equation}
	\mathop{\mathrm{lip}}f(x):=\begin{cases}
		\limsup\limits_{y\to x}\frac{\lvert f(y)-f(x)\rvert}{d(x,y)},\quad & \text{if } x\text{ is not isolated},\\
		0,\quad & \text{otherwise}.
	\end{cases}
\end{equation}
For each $f\in L^{2}(X,m)$, the Cheeger energy is defined by
\begin{equation}
	\mathop{\mathrm{Ch}}(f):=\inf\left\{\liminf_{n\to\infty}\int_{X}(\mathop{\mathrm{lip}}f)^{2}\mathrm{d}m: \lvert f-f_{n}\rvert_{L^{2}(X,m)}\to 0\text{ and } f_{n}\in\mathop{\mathrm{Lip}}(X)\right\}.
\end{equation}
The Sobolev space $W^{1,2}(X,d,m)$ is defined to be the collection of all $f\in L^{2}(X,m)$ with $\mathop{\mathrm{Ch}}(f)<\infty$ and the Sobolev norm is defined by
\begin{equation}
	\lvert f\rvert_{W^{1,2}(X,d,m)}:=\left(\int_{X}f^{2}\mathrm{d}m+\mathop{\mathrm{Ch}}(f)\right)^{\frac{1}{2}}.
\end{equation}
Recall that $RCD(K,N)$-spaces satisfy the doubling property \cite[Theorem 2.3]{Sturm_2006_2} and the local Poincar\'{e} inequality \cite{Rajala_2012}. In particular, the theory developed in \cite{Cheeger_1999} can be applied to $RCD(K,N)$-spaces. Cheeger energy is a convex and lower semi-continuous functional on $L^{2}(X,m)$ and for each $f\in W^{1,2}(X,d,m)$ there exists a unique non-negative function in $L^{2}(X,m)$, denoted by $\lvert \nabla f\rvert$, such that
\begin{equation}
	\mathop{\mathrm{Ch}}(f)=\int_{X}\lvert \nabla f\rvert^{2}\mathrm{d}m.
\end{equation}
$\lvert \nabla f \rvert$ is called the minimal weak upper gradient of $f$. Moreover, it holds that $\mathop{\mathrm{lip}}f=\lvert \nabla f\rvert$ for each $f\in\mathop{\mathrm{Lip}}(X)$.

On $RCD(K,N)$-spaces, Cheeger energy is assumed to be a quadratic form. The inner product $\nabla f \cdot \nabla g$ is well defined for each pair of $f,g\in W^{1,2}(X,d,m)$ and is a symmetric bilinear form \cite{Ambrosio-Gigli-Savare_2014}. Moreover, it holds pointwise that $\lvert \nabla (f+g)\rvert^{2}=\lvert \nabla f\rvert^{2}+ 2\nabla f\cdot\nabla g + \lvert \nabla g\rvert^{2}$ for each pair of $f,g\in W^{1,2}(X,d,m)$.

Given an open domain $\Omega\subset X$, we denote by $\mathop{\mathrm{Lip}}_{0}(\Omega)$ the collection of Lipschitz functions with compact support contained in $\Omega$. The Sobolev space $W^{1,2}_{0}(\Omega)$ is defined to be the closure of $\mathop{\mathrm{Lip}}_{0}(\Omega)$ with respect to the Sobolev norm $\lvert\cdot\rvert_{W^{1,2}(X,d,m)}$. The local Sobolev space $W^{1,2}_{\mathrm{loc}}(\Omega)$ is defined to be the collection of all functions $f$ satisfying that $fg\in W^{1,2}(X,d,m)$ for each $g\in\mathop{\mathrm{Lip}}_{0}(\Omega)$. The Sobolev space $W^{1,2}(\Omega)$ is defined to be the collection of all functions $f\in W^{1,2}_{\mathrm{loc}}(X,d,m)$ such that $f,\lvert \nabla f\rvert\in L^{2}(\Omega,m)$. Recall that the minimal weak upper gradient has locality, i.e. $\lvert\nabla f\rvert=\lvert\nabla g\rvert$ $m$-a.e. on the set $\{f=g\}$ for all $f,g\in W^{1,2}_{\mathrm{loc}}(X,d,m)$.

We recall the notion of measure-valued Laplacian \cite{Gigli_2015}.
\begin{definition}\label{definition_1_preliminaries}
	Let $\Omega\subset X$ be an open domain. For each $f\in W^{1,2}(\Omega)$, the distributional Laplacian of $f$ on $\Omega$, denoted by $\mathbf{\Delta} f$, is the linear functional defined by
	\begin{equation}
		\mathbf{\Delta}f(\phi):=-\int_{\Omega}\nabla f\cdot\nabla \phi\mathrm{d}m,\quad\forall \phi\in\mathop{\mathrm{Lip}_{0}}(\Omega).
	\end{equation}
	We say that $f$ has a measure-valued Laplacian if there exists a Radon measure $\mu$ such that
	\begin{equation}
		\mathbf{\Delta}f(\phi)=\int_{\Omega}\phi\mathrm{d}\mu.
	\end{equation}
	In this case, $\mu$ is unique and we denote by $\mathbf{\Delta}f:=\mu$. In the case that $\mathbf{\Delta}u$ is absolutely continuous with respect to $m$, we denote by $\Delta u$ the density of $\mathbf{\Delta}u$ with respect to $m$, i.e. $\mathbf{\Delta} u=(\Delta u) m$.
\end{definition}
On $RCD(K,N)$-spaces, the following Laplacian comparison theorem for distance functions holds \cite[Corollary 5.15]{Gigli_2015}.
\begin{theorem}\label{theorem_1_preliminaries}
	Let $(X,d,m)$ be an $RCD(K,N)$-space with $K\in\mathbb{R}$ and $N\in(1,\infty)$. Given a point $x_{0}\in X$, the distance function $d_{x_{0}}(x):=d(x,x_{0})$ has a measure-valued Laplacian on $X\setminus\{x_{0}\}$. Moreover, it holds that
	\begin{equation}
		\mathbf{\Delta}d_{x_{0}} \le \frac{\mathfrak{c}_{K,N}(d_{x_{0}})}{d_{x_{0}}}m\quad\text{on } X\setminus\{x_{0}\},
	\end{equation}
	and that
	\begin{equation}
		\mathbf{\Delta} \frac{d_{x_{0}}^{2}}{2} \le \left(\mathfrak{c}_{K,N}(d_{x_{0}})+1\right) m\quad\text{on } X,
	\end{equation}
	where
	\begin{equation}
		\kappa:=\frac{K}{N-1},\qquad \mathfrak{c}_{K,N}(t):=\begin{cases}t(N-1)
			\sqrt{\kappa}\cot(t\sqrt{\kappa}),\quad & \text{if } \kappa>0,\\
			N-1,\quad & \text{if } \kappa=0,\\
			t(N-1)\sqrt{-\kappa}\coth(t\sqrt{-\kappa}),\quad & \text{if } \kappa<0.
		\end{cases}
	\end{equation}
\end{theorem}
\begin{remark}\label{remark_1_preliminaries}
	If $K<0$, then it holds that
	\begin{equation}
		\mathbf{\Delta}\phi_{K,N}(d_{x_{0}})\le 0\quad \text{on }X\setminus\{x_{0}\},
	\end{equation}
	where
	\begin{equation}
	    \phi_{K,N}(t):=-\int_{t}^{1}\left(\frac{\sinh(s\sqrt{-\kappa})}{\sqrt{-\kappa}}\right)^{1-N}\mathrm{d}s.
	\end{equation}
	Indeed, by the chain rule, we have that
	\begin{equation}
		\mathbf{\Delta}\phi_{K,N}(d_{x_{0}})=\phi_{K,N}^{\prime\prime}(d_{x_{0}})\lvert\nabla d_{x_{0}}\rvert^{2} + \phi_{K,N}^{\prime}(d_{x_{0}})\mathbf{\Delta}d_{x_{0}}.
	\end{equation}
	A direct computation gives that
	\begin{equation}
		\begin{split}
			&\phi_{K,N}^{\prime}(t)=\left(\frac{\sinh(t\sqrt{-\kappa})}{\sqrt{-\kappa}}\right)^{1-N},\\
			&\phi_{K,N}^{\prime\prime}(t)=-(N-1)\sqrt{-\kappa}\coth(t\sqrt{-\kappa})\phi_{K,N}^{\prime}(t).
		\end{split}
	\end{equation}
	Note that $\lvert\nabla d_{x_{0}}\rvert=\mathop{\mathrm{lip}}d_{x_{0}}=1$. Therefore, we get that
	\begin{equation}
		\mathbf{\Delta}\phi_{K,N}(d_{x_{0}}) \le 0.
	\end{equation}
\end{remark}
Given $f\in L^{\infty}(\Omega,m)$ and $g\in W^{1,2}(\Omega)$, by minimizing the functional $\int_{\Omega}(\lvert\nabla u\rvert^{2} + 2fu)\mathrm{d}m$, the following Poisson equation if solvable:
\begin{equation}\label{equation_1_preliminaries}
	\begin{cases}
		\Delta u=f,\quad \text{in } \Omega,\\
		u-g\in W^{1,2}_{0}(\Omega).
	\end{cases}
\end{equation}
See also \cite[Theorem 7.12]{Cheeger_1999} for harmonic functions. Moreover, $u$ has Lipschitz regularity \cite{Jiang-Koskela-Yang_2014} and satisfies the comparison principle and the strong maximum principle \cite[Theorem 9.13, Theorem 9.39]{Bjorn-Bjorn_2011}. See also \cite[Theorem 7.17]{Cheeger_1999} for harmonic functions. Note that the comparison principle implies the uniqueness of solutions to the Poisson equations. We summarize these properties as the following theorem.
\begin{theorem}\label{theorem_2_preliminaries}
	Let $(X,d,m)$ be an $RCD(K,N)$-space with $K\in\mathbb{R}$ and $N\in(1,\infty)$ and let $\Omega\subset X$ be a bounded, open domain. The following statements hold.
	\begin{itemize}
		\item[(1)] Given $f\in L^{\infty}(\Omega)$ and $g\in W^{1,2}(\Omega)$, there exists a unique solution to the Poisson equation \eqref{equation_1_preliminaries}.
		\item[(2)] If $\Delta u \in L^{\infty}(\Omega)$, then $u$ is locally Lipschitz continuous in $\Omega$.
		\item[(3)] If $u\in W^{1,2}_{0}(\Omega)$ is subharmonic, i.e. $\Delta u\ge 0$, then $u\le 0$ in $\Omega$.
		\item[(4)] If $u$ is subharmonic in $\Omega$ and attains its maximum in $\Omega$, then $u$ is constant in $\Omega$.
	\end{itemize}
\end{theorem}
Yau's gradient estimate also holds on $RCD(K,N)$-spaces \cite{Zhang-Zhu_2016}.
\begin{theorem}\label{theorem_3_preliminaries}
	Let $(X,d,m)$ be an $RCD(K,N)$-space with $K\in\mathbb{R}$ and $N\in(1,\infty)$ and let $u$ be a positive harmonic function on $B_{2r}(x_{0})\subset X$. Then the following Yau's gradient estimate holds
	\begin{equation}
		\sup_{B_{r}(x_{0})}\frac{\lvert\nabla u\rvert}{u}\le \frac{C(1+\sqrt{\lvert K\rvert}r)}{r},
	\end{equation}
	where the constant $C>0$ depends only on $N$.
\end{theorem}
\begin{remark}\label{remark_2_preliminaries}
	Note that Yau's gradient estimate implies Harnack inequality. More precisely, for a positive harmonic function $u$ on $B_{2r}(x_{0})$ it holds that
	\begin{equation}
		u(x)\le e^{C_{N}(1+\sqrt{\lvert K\rvert}r)}u(y),\quad\forall x,y\in B_{r}(x_{0}).
	\end{equation}
	Indeed, let $\gamma$ be a geodesic connected $x$ and $y$. Then
	\begin{equation}
		\ln u(x)-\ln u(y)=\int_{\gamma} \frac{\mathrm{d}}{\mathrm{d}t}\ln u(\gamma_{t})\le \left(\sup_{B_{r}(x_{0})}\frac{\lvert\nabla u\rvert}{u}\right) d(x,y)\le C_{N}(1+\sqrt{\lvert K\rvert}r).
	\end{equation}
	Equivalently,
	\begin{equation}
		u(x)\le e^{C_{N}(1+\sqrt{\lvert K\rvert}r)}u(y).
	\end{equation}
\end{remark}
On $RCD(K,N)$-spaces, the existence of good cut-off functions has been established \cite[Lemma 3.1]{Mondino-Naber_2019}.
\begin{lemma}\label{lemma_1_preliminaries}
	Let $(X,d,m)$ be an $RCD(K,N)$-space with $K\in\mathbb{R}$ and $N\in(1,\infty)$. Then for every $x\in X$ and every $0<r<R$ there exists a Lipschitz function $\phi:X\to\mathbb{R}$ satisfying that
	\begin{itemize}
		\item[(1)] $0\le\phi\le 1$ on $X$, $\phi=1$ on $B_{r}(x)$ and $\phi=0$ on $X\setminus B_{2r}(x)$.
		\item[(2)] $r^{2}\lvert\Delta\phi\rvert+r\lvert\nabla\phi\rvert\le C$ where the constant $C>0$ depends only on $K,N,R$.
	\end{itemize}
\end{lemma}
We recall the notions of functions of bounded variation and perimeter measure and the coarea formula \cite{Miranda_2003}. A function $f\in L^{1}(X,m)$ is said to be of bounded variation, denoted by $f\in \mathop{\mathrm{BV}}(X,d,m)$ if there exists a sequence of locally Lipschitz continuous functions $f_{k}$ converging to $f$ in $L^{1}(X,m)$ such that
\begin{equation}
	\limsup_{k\to\infty} \int_{X}\lvert\nabla f_{k}\rvert\mathrm{d}m<\infty.
\end{equation}
For each open $A\subset X$, we define
\begin{equation}
	\lvert Df\rvert(A):=\inf \liminf_{k\to\infty} \int_{A}\lvert\nabla f_{k}\rvert\mathrm{d}m,
\end{equation}
where the infimum is taken over all sequences of locally Lipschitz continuous functions $f_{k}$ converging to $f$ in $L^{1}(X,m)$. It can be proved that $\lvert Df\rvert$ can be extended to a Borel measure on $(X,d)$. We call this Borel measure the total variation of $f$ and still denote it by $\lvert Df\rvert$.

Given a Borel set $E\subset X$ and an open set $A\subset X$, the perimeter of $E$ on $A$ is defined to be
\begin{equation}
	\mathop{\mathrm{Per}}(E;A):=\inf\liminf_{k\to\infty} \int_{A} \lvert \nabla f_{k}\rvert\mathrm{d}m,
\end{equation}
where the infimum is taken over all sequences of locally Lipschitz continuous functions $f_{k}:A\to\mathbb{R}$ such that $f_{k}$ converges to $\chi_{E}$ in $L^{1}_{\mathrm{loc}}(A,m)$. If $\mathop{\mathrm{Per}}(E;X)<\infty$, then we say that $E$ has finite perimeter. In this case, it can be proved that $\mathop{\mathrm{Per}}(E;\cdot)$ can be extended to a Borel measure on $(X,d)$. Moreover, the perimeter measure is lower semi-continuous by the definition. That is, 
\begin{equation}
	\mathop{\mathrm{Per}}(E;A) \le \liminf_{k\to\infty} \mathop{\mathrm{Per}}(E_{k};A)
\end{equation}
whenever $\chi_{E_{k}}$ converges to $\chi_{E}$ in $L^{1}(A,m)$. We say that $E$ has locally finite perimeter if for each $x\in X$ there exists an $r>0$ such that $\mathop{\mathrm{Per}}\left(E;B_{r}(x)\right)<\infty$.

We will need the following coarea formula \cite[Remark 4.3]{Miranda_2003}.
\begin{theorem}\label{theorem_5_preliminaries}
	Let $(X,d,m)$ be an $RCD(K,N)$-space with $K\in\mathbb{R}$ and $N\in(1,\infty)$ and let $v\in\mathop{\mathrm{BV}}(X,d,m)$ be continuous and non-negative. Then it holds that
	\begin{equation}
		\int_{\{s\le v<t\}} f\mathrm{d}\lvert Dv\rvert = \int_{s}^{t} \int_{X} f\mathrm{d}\mathop{\mathrm{Per}}(\{v>r\};\cdot) \mathrm{d}r
	\end{equation}
	for all Borel functions $f:X\to[0,\infty]$ and all $0\le s<t<\infty$.
\end{theorem}

\subsection{Pointed measured Gromov-Hausdorff convergence and stability results}
We recall the notion of pointed measured Gromov-Hausdorff convergence \cite{Gigli-Mondino-Savare_2015}. We always assume that $\mathop{\mathrm{supp}}m=X$ for each metric measure space $(X,d,m)$.
\begin{definition}\label{definition_2_preliminaries}
	We say that a sequence of pointed metric measure spaces $(X_{k},d_{k},m_{k},x_{k})$ converges to $(X_{\infty},d_{\infty},m_{\infty},x_{\infty})$ in the pointed measured Gromov-Hausdorff topology (p-mGH for short) if there exist a separable metric space $(Z,d)$ and isometric embeddings $\phi_{k}:X_{k}\to Z$, $k\in\mathbb{N}\cup\{\infty\}$ satisfying the following properties.
	\begin{itemize}
		\item[(1)] For every $\varepsilon>0$ and $R>0$ there exists a $k_{0}$ such that
		\begin{equation}
			\phi_{\infty}\left(B_{R}^{X_{\infty}}(x_{\infty})\right) \subset \mathcal{N}_{\varepsilon}\left(\phi_{k}\big(B_{R}^{X_{k}}(x_{k})\big)\right), \quad \text{and} \quad \phi_{k}\left(B_{R}^{X_{k}}(x_{k})\right) \subset \mathcal{N}_{\varepsilon} \left(\phi_{\infty}\big(B_{R}^{X_{\infty}}(x_{\infty})\big)\right),
		\end{equation}
		for all $k\ge k_{0}$, where $\mathcal{N}_{\varepsilon}(A)$ is the $\varepsilon$-neighborhood of $A$ for each $A\subset Z$, i.e. $\mathcal{N}_{\varepsilon}(A):=\{z\in Z: d(z,A)<\varepsilon\}$.
		\item[(2)] For each bounded continuous function $\varphi:Z\to\mathbb{R}$ with bounded support, it holds that
		\begin{equation}
			\lim_{k\to\infty}\int_{Z}\varphi\mathrm{d}(\phi_{k})_{\#}m_{k} = \int_{Z}\varphi \mathrm{d}(\phi_{\infty})_{\#}m_{\infty}.
		\end{equation}
	\end{itemize}
\end{definition}
\begin{remark}\label{remark_3_preliminaries}
	By choosing a $\varphi:Z\to[0,1]$ satisfying that $\varphi=1$ on $B_{R+\varepsilon}$ and $\varphi=0$ on $Z\setminus B_{2R}$, we get that $\limsup\limits_{k\to\infty} m_{k}\big(B_{R}^{X_{k}}(x_{k})\big)\le m_{\infty}\big(B_{2R}^{X_{\infty}}(x_{\infty})\big)$.
\end{remark}
Recall that a pointed metric measure space $(X,d,m,x)$ is said to be normalized if it holds that
\begin{equation}
	\int_{B_{1}(x)} \big(1-d(\cdot,x)\big) \mathrm{d}m = 1.
\end{equation}

By \cite[Theorem 3.1]{Sturm_2006_2} and \cite{Gigli-Mondino-Savare_2015}, we have the following compactness result.
\begin{theorem}\label{theorem_4_preliminaries}
	Let $K\in\mathbb{R}$ and $N\in(1,\infty)$. Then the class of normalized pointed $RCD(K,N)$-spaces is compact with respect to p-mGH convergence.
\end{theorem}

We recall the notions of convergence of functions defined on varying metric measure spaces and two stability results that will be used later \cite{Ambrosio-Honda_2017,Ambrosio-Honda_2018}. In the rest of this subsection,  $(X_{k},d_{k},m_{k},x_{k})$ is a sequence of pointed $RCD(K,N)$-spaces with $K\in\mathbb{R},N\in(1,\infty)$ converging to $(X_{\infty},d_{\infty},m_{\infty},x_{\infty})$ in the p-mGH topology and $(Z,d)$ is the metric space as in Definition \ref{definition_2_preliminaries}.
\begin{definition}\label{definition_3_preliminaries}
	We say that $f_{k}\in L^{2}(X_{k},m_{k})$ $L^{2}$-weakly converges to $f_{\infty}\in L^{2}(X_{\infty},m_{\infty})$ if $\sup\limits_{k} \lvert f_{k}\rvert_{L^{2}}<\infty$ and $f_{k}m_{k}\to f_{\infty}m_{\infty}$ in duality to bounded continuous functions on $(Z,d)$ with bounded support. We say that $f_{k}$ $L^{2}$-strongly converges to $f_{\infty}$  if $f_{k}$ $L^{2}$-weakly converges to $f_{\infty}$ with $\lim\limits_{k\to\infty} \lvert f_{k}\rvert_{L^{2}}= \lvert f_{\infty }\rvert_{L^{2}}$.
\end{definition}
\begin{definition}\label{definition_4_preliminaries}
	We say that $f_{k}\in W^{1,2}(X_{k},d_{k},m_{k})$ $W^{1,2}$-weakly converges to $f_{\infty}\in W^{1,2}(X_{\infty},d_{\infty},m_{\infty})$ if $f_{k}$ $L^{2}$-weakly converges to $f_{\infty}$ and $\sup\limits_{k}\lvert \nabla f_{k}\rvert_{L^{2}}<\infty$. We say that $f_{k}$ $W^{1,2}$-strongly converges to $f_{\infty}$ if $f_{k}$ $L^{2}$-strongly converges to $f_{\infty}$ and $\lim\limits_{k\to\infty} \lvert \nabla f_{k}\rvert_{L^{2}}=\lvert\nabla f_{\infty}\rvert_{L^{2}}$.
\end{definition}
\begin{lemma}[\text{\cite[Lemma 2.10]{Ambrosio-Honda_2018}}]\label{lemma_2_preliminaries}
	For each $f_{\infty}\in\mathop{\mathrm{Lip}_{0}}\big(B_{R}^{\infty}(x_{\infty})\big)$, there exists a sequence of $f_{k}\in\mathop{\mathrm{Lip}_{0}}\big(B_{R}^{k}(x_{k})\big)$ $W^{1,2}$-strongly converging to $f_{\infty}$ with $\sup\limits_{k} \lvert \nabla f_{k}\rvert_{L^{\infty}}<\infty$.
\end{lemma}
\begin{lemma}[\text{\cite[Theorem 1.5.7]{Ambrosio-Honda_2017}}]\label{lemma_3_preliminaries}
	Let $f_{k}\in W^{1,2}(X_{k},d_{k},m_{k})$ $W^{1,2}$-strongly converge to $f_{\infty}\in W^{1,2}(X_{\infty},d_{\infty},m_{\infty})$ and let $g_{k}\in W^{1,2}(X_{k},d_{k},m_{k})$ $W^{1,2}$-weakly converge to $g_{\infty}\in W^{1,2}(X_{\infty},d_{\infty},m_{\infty})$. If $\sup\limits_{k} \lvert \nabla f_{k}\cdot\nabla g_{k}\rvert_{L^{2}}<\infty$, then $\nabla f_{k}\cdot\nabla g_{k}$ $L^{2}$-weakly converges to $\nabla f_{\infty}\cdot\nabla g_{\infty}$.
\end{lemma}

\section{Existence}
Since the energy functional \eqref{equation_2_Introduction} is lower semi-continuous and the set $\mathscr{A}_{g}$ is weakly convex, the existence of a minimizer follows from the classical theory of the calculus of variations. For completeness, we include a proof in this section.
\begin{lemma}\label{lemma_Existence}
	Let $(X,d,m)$ be an $RCD(K,N)$-space with $K\in\mathbb{R}$ and $N\in(1,\infty)$ and let $\Omega \subset X$ be a bounded and open domain. Given $g\in W^{1,2}(\Omega)$, we define $\mathscr{A}_{g}:=\{v\in W^{1,2}(\Omega):v-g\in W^{1,2}_{0}(\Omega)\}$. Then there exists a minimizer of $E$ over $\mathscr{A}_{g}$. More precisely, there exists a $u\in\mathscr{A}_{g}$ such that
	\begin{equation}
		E(u)\le E(v),\quad \forall v\in\mathscr{A}_{g}.
	\end{equation}
\end{lemma}
\begin{proof}
	Note that
	\begin{equation}
		0\le\inf_{v\in\mathscr{A}_{g}}E(v)\le E(g) \le \int_{\Omega}\lvert \nabla g\rvert^{2}\mathrm{d}m + (\lambda_{+}+\lambda_{-})m(\Omega)<\infty.
	\end{equation}
	Let $\{v_{k}\}\subset\mathscr{A}_{g}$ be a minimizing sequence, i.e.
	\begin{equation}
		\lim_{k\to\infty}E(v_{k})=\inf_{v\in\mathscr{A}_{g}}E(v).
	\end{equation}
	Then $\sup\limits_{k}\int_{\Omega}\lvert \nabla v_{k}\rvert^{2}\mathrm{d}m<\infty$. By the local Poincar\'{e} inequality \cite{Rajala_2012}, we get that $\{v_{k}\}$ is a bounded sequence in $W^{1,2}(\Omega)$. By the weak compactness of $W^{1,2}(\Omega)$, $v_{k}$ converges to some $u$ weakly in $W^{1,2}(\Omega)$ and $m$-a.e. on $\Omega$, up to a subsequence. Moreover, $u\in\mathscr{A}_{g}$ since $\mathscr{A}_{g}$ is weakly closed.
	
	By the lower semi-continuity of Cheeger energy, we have that
	\begin{equation}
		\int_{\Omega}\lvert\nabla u\rvert^{2}\mathrm{d}m\le\liminf_{k\to\infty}\int_{\Omega} \lvert\nabla v_{k}\rvert^{2}\mathrm{d}m.
	\end{equation}
	By Fatou's lemma, we have that
	\begin{equation}
		\int_{\Omega} \lambda_{+}\chi_{\{u>0\}}+\lambda_{-}\chi_{\{u<0\}}\mathrm{d}m\le\liminf_{k\to\infty}\int_{\Omega}\lambda_{+}\chi_{\{v_{k}>0\}}+\lambda_{-}\chi_{\{v_{k}<0\}}\mathrm{d}m.
	\end{equation}
	Therefore, we get that
	\begin{equation}
		E(u)\le\liminf_{k\to\infty}E(v_{k})=\inf_{v\in\mathscr{A}_{g}} E(v).
	\end{equation}
\end{proof}

\section{H\"{o}lder continuity}
In this section, we prove that minimizers are H\"{o}lder continuous. The proof is based on the Campanato theorem \cite{Gorka_2009}.
\begin{lemma}\label{lemma_1_Holder_continuouity}
	Let $u$ be a minimizer of $E$. Then $\lvert u\rvert$ is subharmonic in $\Omega$, i.e. $\mathbf{\Delta}\lvert u\rvert\ge 0$ in $\Omega$.
\end{lemma}
\begin{proof}
	For each $0\le\phi\in\mathop{\mathrm{Lip}_{0}}(\Omega)$ and each $t>0$, we define
	\begin{equation}
		u_{t}(x):=\begin{cases}
			\max\{u(x)-t\phi(x),0\},\quad& \text{if } u(x)>0,\\
			\min\{u(x)+t\phi(x),0\},\quad& \text{if } u(x)\le 0.
		\end{cases}
	\end{equation}
	Then we have that $u_{t}\in\mathscr{A}_{g}$ since $\phi\in\mathop{\mathrm{Lip}_{0}}(\Omega)$. Hence we have that
	\begin{equation}
		E(u_{t})\ge E(u).
	\end{equation}
	On the other hand, we have that $\chi_{\{u_{t}>0\}}\le\chi_{\{u>0\}}$ and that $\chi_{\{u_{t}<0\}}\le\chi_{\{u<0\}}$ by the definition of $u_{t}$. Therefore, we get that
	\begin{equation}
		\begin{split}
			0 \le & \int_{\Omega}\lvert\nabla u_{t}\rvert^{2}\mathrm{d}m - \int_{\Omega}\lvert\nabla u\rvert^{2}\mathrm{d}m\\
			= & \int_{\Omega} \lvert\nabla (u-t\phi)\rvert^{2}\chi_{\{u>0\}} + \lvert\nabla(u+t\phi)\rvert^{2}\chi_{\{u<0\}}\mathrm{d}m-\int_{\Omega}\lvert\nabla u\rvert^{2}\mathrm{d}m\\
			= & -2t\int_{\Omega} \nabla\phi\cdot\nabla\lvert u\rvert\mathrm{d}m + t^{2}\int_{\Omega}\lvert\nabla\phi\rvert^{2}\mathrm{d}m.
		\end{split}
	\end{equation}
	By letting $t\to 0$, we get that
	\begin{equation}
		-\int_{\Omega}\nabla\phi\cdot\nabla\lvert u\rvert\mathrm{d}m\ge 0.
	\end{equation}
	By Riesz representation theorem, we get that $\mathbf{\Delta}\lvert u\rvert\ge 0$ in $\Omega$.
\end{proof}
\begin{remark}
	By taking $u_{t}:=ue^{-t\phi}$, it can be proved that $\Delta u=2\lvert\nabla u\rvert^{2}$.
\end{remark}
\begin{corollary}\label{corollary_1_Holder_continuity}
	Let $u$ be a minimizer of $E$. Then $u$ is locally bounded in $\Omega$. More precisely, for each $R>0$ there exists a constant $C>0$, depending only on $K,N,R$, such that
	\begin{equation}
		\sup_{B_{r}(x)} \lvert u\rvert \le C\left(\fint_{B_{2r}(x)}\lvert u\rvert^{2}\mathrm{d}m\right)^{\frac{1}{2}},\quad\forall B_{2R}(x)\subset\Omega,\forall r\in(0,R).
	\end{equation}
\end{corollary}
\begin{proof}
	It is a direct consequence of Lemma \ref{lemma_1_Holder_continuouity} and \cite[Theorem 4.3]{Kinnunen-Shanmugalingam_2001}.
\end{proof}
The main result of this section is the following theorem.
\begin{theorem}\label{theorem_1_Holder_continuity}
	Let $u$ be a minimizer of $E$ and let $0<\lambda_{-}<\lambda_{+}\le\lambda_{0}$. Then $u$ is locally H\"{o}lder continuous in $\Omega$. More precisely, for each ball $B_{2r}(x_{0})\subset \Omega$ with $r^{1+\frac{2}{N}}<\frac{r}{2}$ there exists a constant $C>0$, depending only on $K,N,\lambda_{0}$ and $\lvert u\rvert_{L^{\infty}(B_{2r}(x_{0}))}$, such that
	\begin{equation}
		\lvert u(x)-u(y)\rvert\le Cd(x,y)^{\frac{2}{N+2}},\quad\forall x,y\in B_{r}(x_{0}).
	\end{equation}
\end{theorem}
\begin{proof}
	Let $B_{2r}(x_{0})\subset\Omega$ with $r^{1+\frac{2}{N}}<\frac{r}{2}$ and let $v$ be the harmonic function on $B_{r}(x_{0})$ with boundary value $v-u\in W^{1,2}_{0}(B_{r}(x_{0}))$. By the maximum principle, we have that
	\begin{equation}
		\lvert v \rvert_{L^{\infty}(B_{r}(x_{0}))}\le \sup_{\partial B_{r}(x_{0})}\lvert v\rvert=\sup_{\partial B_{r}(x_{0})}\lvert u\rvert\le C_{K,N}\lvert u\rvert_{L^{2}(B_{2r}(x_{0}))}=:M.
	\end{equation}
	By Yau's gradient estimate, we get that
	\begin{equation}
		\lvert\nabla v\rvert \le \frac{C_{K,N}M}{r},\quad\text{on } B_{\frac{r}{2}}(x_{0}).
	\end{equation}
	We extend $v$ to $\Omega$ by letting $v=u$ on $\Omega\setminus B_{r}(x_{0})$. Then $v\in\mathscr{A}_{g}$ and hence $E(u)\le E(v)$. In particular, we get that
	\begin{equation}
		\begin{split}
			\int_{B_{r}(x_{0})} \lvert \nabla u\rvert^{2}-\lvert\nabla v\rvert^{2}\mathrm{d}m \le &\int_{B_{r}(x_{0})} \lambda_{+}(\chi_{\{v>0\}}-\chi_{\{u>0\}})+\lambda_{-}(\chi_{\{v<0\}}-\chi_{\{u<0\}})\mathrm{d}m\\
			\le & \lambda_{0}m(B_{r}(x_{0})).
		\end{split}
	\end{equation}
	On the other hand, since $v$ is harmonic in $B_{r}(x_{0})$ with $u-v\in W^{1,2}_{0}(B_{r}(x_{0}))$, we have that
	\begin{equation}
		\int_{B_{r}(x_{0})}\lvert\nabla u\rvert^{2}-\lvert\nabla v\rvert^{2}\mathrm{d}m=\int_{B_{r}(x_{0})}\nabla(u+v)\cdot\nabla(u-v)\mathrm{d}m=\int_{B_{r}(x_{0})} \lvert\nabla(u-v)\rvert^{2}\mathrm{d}m.
	\end{equation}
	For $\rho\in(0,\frac{r}{2})$, we estimate $\int_{B_{\rho}(x_{0})}\lvert\nabla u\rvert^{2}\mathrm{d}m$ as below.
	\begin{equation}\label{equation_1_Holder_continuity}
		\begin{split}
			\int_{B_{\rho}(x_{0})}\lvert\nabla u\rvert^{2}\mathrm{d}m \le & 2\int_{B_{\rho}(x_{0})} \lvert\nabla(u-v)\rvert^{2} + \lvert\nabla v\rvert^{2}\mathrm{d}m\\
			\le & 2\int_{B_{r}(x_{0})} \lvert \nabla u\rvert^{2}-\lvert\nabla v\rvert^{2}\mathrm{d}m + 2\int_{B_{\rho}(x_{0})} \frac{C_{K,N}M^{2}}{r^{2}}\mathrm{d}m\\
			\le & \lambda_{0}m(B_{r}(x_{0})) + C_{K,N}\frac{M^{2}}{r^{2}}m(B_{\rho}(x_{0}))\\
			\le & C_{K,N,\lambda_{0}}m(B_{\rho}(x_{0}))\left(\frac{r^{N}}{\rho^{N}}+\frac{M^{2}}{r^{2}}\right),
		\end{split}
	\end{equation}
	where we use the doubling property in the last inequality. By choosing $\rho=r^{1+\frac{2}{N}}$, we get that
	\begin{equation}
		\rho^{2}\fint_{B_{\rho}(x_{0})}\lvert\nabla u\rvert^{2}\mathrm{d}m \le C_{K,N,\lambda_{0}}\left(\frac{r^{N}}{\rho^{N-2}}+\frac{M^{2}\rho^{2}}{r^{2}}\right)=C_{K,N,\lambda_{0}}(1+M^{2})\rho^{\frac{4}{N+2}}.
	\end{equation}
	By the Campanato theorem \cite[Theorem 3.3]{Gorka_2009}, we get that
	\begin{equation}
		\lvert u(x)-u(y)\rvert\le C_{K,N,\lambda_{0},M}d(x,y)^{\frac{2}{N+2}},\quad\forall x,y\in B_{r}(x_{0}).
	\end{equation}
\end{proof}
\begin{corollary}\label{corollary_2_Holder_continuity}
	Let $u$ be minimizer of $E$. Then the sets $\{u_{\pm}>0\}$ are open subsets of $\Omega$. In particular, $u$ is harmonic (hence locally Lipschitz continuous) in $\{u_{\pm}>0\}$.
\end{corollary}
\begin{proof}
	The openness of $\{u_{\pm}>0\}$ follows from the continuity of $u$. For each $B_{2r}\subset\{u>0\}$ and each $\phi\in\mathop{\mathrm{Lip}_{0}}(B_{r})$, we define $u_{t}:=u+t\phi$, $t\in\mathbb{R}$. Then $u_{t}\in\mathscr{A}_{g}$ and hence $E(u)\le E(u_{t})$. Note that $\inf_{B_{r}}u>0$. Therefore, for sufficiently small $t$, we have that $u_{t}>0$ in $B_{r}$, which implies that
	\begin{equation}
		\int_{B_{r}}\lvert\nabla
		 u\rvert^{2}\mathrm{d}m\le\int_{B_{r}}\lvert\nabla u_{t}\rvert^{2}\mathrm{d}m.
	\end{equation}
	By letting $t\to 0$, we get that
	\begin{equation}
		-\int_{B_{r}}\nabla u\cdot\nabla\phi\mathrm{d}m=0.
	\end{equation}
	That is, $u$ is harmonic in $B_{r}$ and hence in $\{u>0\}$. A similar argument gives that $u$ is harmonic in $\{u<0\}$.
\end{proof}

\section{Lipschitz continuity}
In this section, we investigate the Lipschitz continuity of minimizers. We employ the method from \cite{Danielli-Petrosyan_2005}, which originally dealt with $p$-Laplacian operator on $\mathbb{R}^{n}$. Compared to \cite{Danielli-Petrosyan_2005}, we need to deal with Sobolev functions defined on varying metric measure spaces. Consequently, we rely on the stability results in \cite{Ambrosio-Honda_2017}.
\begin{lemma}\label{lemma_1_Lipschitz_continuity}
	Let $v$ be a non-negative continuous function in $B_{2}(x_{0})$ with $x_{0}\in\partial\{v>0\}$ and be harmonic in $\{v>0\}$. Let $M:=\sup\{v(x):x\in B_{\frac{1}{20}}(x_{0})\}$. Then there exist $y\in B_{1}(x_{0})$ and $\eta\in B_{1}(x_{0})\cap\partial\{v>0\}$ such that
	\begin{itemize}
		\item[(1)] $v(y)\ge M$.
		\item[(2)] $d(y,\eta)=d(y,\{v=0\})=:\rho$.
		\item[(3)] $v\le v(y)$ in $B_{\frac{\rho}{10}}(\eta)$.
		\item[(4)] $\sup\{v(x):x\in B_{\frac{\rho}{20}}(\eta)\}\ge C_{K,N}v(y)$.
	\end{itemize}
\end{lemma}
\begin{proof}
	We define the set
	\begin{equation}
		S:=\left\{x\in \overline{B}_{1}(x_{0}):d(x,\{v=0\})\le\frac{1-d(x,x_{0})}{10}\right\}.
	\end{equation}
	Then $S$ is closed and $B_{\frac{1}{20}}(x_{0})\subset S$. Moreover, there exists a maximum point $y$ of $v$ on $S$, i.e.
	\begin{equation}
		v(y)=\max\{v(x):x\in S\}.
	\end{equation}
	In particular, $v(y)\ge M$. 
	
	Let $\eta\in\{v=0\}$ be such that
	\begin{equation}
		d(y,\eta)=d(y,\{u=0\})=:\rho.
	\end{equation}
	Then we have that
	\begin{equation}
		\rho\le\frac{1-d(y,x_{0})}{10},
	\end{equation}
	which implies that
	\begin{equation}
		B_{\frac{\rho}{10}}(\eta)\subset B_{2\rho}(y)\subset B_{1}(x_{0}).
	\end{equation}
	
	The next step is to show that $B_{\frac{\rho}{10}}(\eta)\subset S$ and thus $v\le M\le v(y)$ in $B_{\frac{\rho}{10}}(\eta)$. Indeed, for $x\in B_{\frac{\rho}{10}}(\eta)$, we have that
	\begin{equation}
		\begin{split}
			d(x,\{v=0\}) \le & d(x,\eta) < \frac{\rho}{10}<\frac{89\rho}{100}\\
			= & \rho-\frac{1}{10}\left(\frac{\rho}{10}+\rho\right)\\
			\le & d(y,\{v=0\})-\frac{d(x,\eta)+d(y,\eta)}{10}\\
			\le & \frac{1-d(y,x_{0})}{10} - \frac{d(x,y)}{10}\\
			\le & \frac{1-d(x,x_{0})}{10}.
		\end{split}
	\end{equation}
	The final step is to show that $\sup\{v(x):x\in B_{\frac{\rho}{20}}(\eta)\}\ge C_{K,N}v(y)$. It suffices to show that $v(x)\ge C_{K,N}v(y)$ for some $x\in\overline{B}_{\frac{\rho}{20}}(x_{0})$. Indeed, let $\gamma:[0,\rho]\to X$ be a geodesic from $\eta$ to $y$ and let $x:=\gamma(\frac{\rho}{20})$. Then $v$ is harmonic in $B_{\rho}(y)\subset\{v>0\}$ and hence, by the Harnack inequality, it holds that $v(x)\ge C_{K,N}v(y)$.
\end{proof}
\begin{lemma}\label{lemma_2_Lipschitz_continuity}
	Let $u$ be a minimizer of $E$ with $0<\lambda_{-}<\lambda_{+}\le\lambda_{0}$. Then it holds that
	\begin{equation}
		-\int_{\Omega}\nabla u\cdot\nabla \phi\mathrm{d}m\le\lvert\nabla\phi\rvert_{L^{2}(\Omega)}\sqrt{\lambda_{0}m(\{\phi\neq 0\})},\quad\forall \phi\in\mathop{\mathrm{Lip}_{0}}(\Omega).
	\end{equation}
\end{lemma}
\begin{proof}
	Since $u+\delta\phi\in\mathscr{A}_{g}$ for all $\delta>0$ and all $\phi\in\mathop{\mathrm{Lip}_{0}}(\Omega)$, we get that
	\begin{equation}
		E(u)\le E(u+\delta\phi).
	\end{equation}
	In particular, we get that
	\begin{equation}
		\begin{split}
			-2\delta\int_{\Omega}\nabla u\cdot\nabla\phi\mathrm{d}m \le & \delta^{2}\int_{\Omega}\lvert\nabla\phi\rvert^{2}\mathrm{d}m + \int_{\Omega}\lambda_{+}(\chi_{\{u+\delta\phi>0\}}-\chi_{\{u>0\}})-\lambda_{-}(\chi_{\{u+\delta\phi<0\}}-\chi_{\{u<0\}})\mathrm{d}m\\
			\le & \delta^{2}\int_{\Omega}\lvert\nabla\phi\rvert^{2}\mathrm{d}m + \lambda_{0}m(\{\phi\neq 0\}).
		\end{split}
	\end{equation}
	The conclusion follows by taking $\delta=\frac{\sqrt{\lambda_{0}m(\{\phi\neq 0\})}}{\lvert\nabla\phi\rvert_{L^{2}(\Omega)}}$.
\end{proof}
We introduce the notion of linear growth.
\begin{definition}\label{definition_1_Lipschitz_continuity}
	Let $v$ be a continuous function on $\Omega$. We say that $v$ has linear growth locally if for each $B_{2R}(x_{0})\subset\Omega$ with $x_{0}\in\{v=0\}$ there exists a constant $L>0$ such that
	\begin{equation}
		\lvert v(x)\rvert \le Ld(x,\{v=0\}),\quad\forall x\in B_{R}(x_{0}).
	\end{equation}
\end{definition}
\begin{lemma}\label{lemma_3_Lipschitz_continuity}
	Let $u$ be a minimizer of $E$. If one of $u_{\pm}$ has linear growth locally, then $u$ has linear growth locally.
\end{lemma}
\begin{proof}
	Without loss of generality, we assume that $u_{-}$ has linear growth locally. Arguing by contradiction, there exist a ball $B_{2R}(x_{0})\subset\Omega$ with $x_{0}\in\{u=0\}$ and a sequence of $x_{k}\in B_{R}(x_{0})\cap\{u>0\}$ such that
	\begin{equation}
		u(x_{k})\ge kd(x_{k},\{u=0\}).
	\end{equation}
	
	\textbf{Step 1.} Let $x_{k}^{0}\in\{u=0\}$ be such that
	\begin{equation}
		d(x_{k},x_{k}^{0})=d(x_{k},\{u=0\})=:r_{k}.
	\end{equation}
	Recall that $u$ is locally bounded. In particular, we have that
	\begin{equation}
		r_{k}\le\frac{1}{k}u(x_{k})\le\frac{1}{k}\sup_{B_{R}(x_{0})}u\to 0,\quad\text{as } k\to\infty.
	\end{equation}
	In particular, we may assume that $r_{k}<10^{-100}R$.
	
	\textbf{Step 2.} By Lemma \ref{lemma_1_Lipschitz_continuity}, for each $x_{k}$ there exist a $y_{k}\in B_{20r_{k}}(x_{k}^{0})$ and an $\eta_{k}\in B_{20r_{k}}(x_{k}^{0})\cap\partial\{u>0\}$ such that
	\begin{equation}\label{equation_1_Lipschitz_continuity}
		\begin{cases}
			u(y_{k})\ge u(x_{k}),\\
			\rho_{k}:=d(y_{k},\{u=0\})=d(y_{k},\eta_{k}),\\
			u\le u(y_{k})\text{ in } B_{\rho_{k}/10}(\eta_{k}),\\
			\sup\{u(x):x\in B_{\rho_{k}/20}(\eta_{k})\}\ge C_{K,N,R}u(y_{k}).
		\end{cases}
	\end{equation}
	We define $d_{k}:=\rho_{k}^{-1}d$, $c_{k}:=\int_{B_{\rho_{k}}(\eta_{k})}\big(1-\rho_{k}^{-1}d(\cdot,\eta_{k})\big)\mathrm{d}m$ and $m_{k}:=c_{k}^{-1}m$. By the compactness, $(X,d_{k},m_{k},\eta_{k})$ p-mGH converges to some pointed metric measure space $(X_{\infty},d_{\infty},m_{\infty},\eta_{\infty})$, up to a subsequence.
	
	We define 
	\begin{equation}
		u_{k}(x):=\frac{u(x)}{u(y_{k})},\quad x\in(X,d_{k},m_{k},\eta_{k}).
	\end{equation}
	Then
	\begin{equation}
		\sup_{B_{1/10}^{k}(\eta_{k})} u_{k} = \sup_{B_{\rho_{k}/10}(\eta_{k})} \frac{u}{u(y_{k})} \le 1,\quad \forall k\in\mathbb{N}.
	\end{equation}
	Here, $B_{r}^{k}(\eta_{k})$ denotes the balls on $(X,d_{k},m_{k},\eta_{k})$.
	
	On the other hand, since $u_{-}$ has linear growth locally, there exists a constant $L>0$ such that
	\begin{equation}
		\inf_{B_{1/10}^{k}(\eta_{k})} u_{k} = \inf_{B_{\rho_{k}/10}(\eta_{k})} \frac{u}{u(y_{k})} \ge -\frac{L\rho_{k}}{10u(y_{k})}.
	\end{equation}
	Recall that $u(y_{k})\ge u(x_{k})\ge kr_{k}\ge \frac{k\rho_{k}}{20}$. Therefore,
	\begin{equation}\label{equation_2_Lipschitz_continuity}
		\inf_{B_{1/10}^{k}(\eta_{k})} u_{k} \ge - \frac{2L}{k}\to 0,\quad\text{as } k\to\infty.
	\end{equation}
	In particular, we may assume that $\inf\limits_{B_{1/10}^{k}(\eta_{k})} u_{k} \ge - 1$ and hence $\lvert u_{k}\rvert\le 1$ in $B_{1/10}^{k}(\eta_{k})$.
	
	\textbf{Step 3.} Note that $\nabla u_{k} = \frac{\rho_{k}\nabla u}{u(y_{k})}$. Therefore,
	\begin{equation}
		\begin{split}
			& \frac{\rho_{k}^{2}}{c_{k}\lvert u(y_{k})\rvert^{2}} \int_{B_{\rho_{k}/10}(\eta_{k})} \lvert \nabla u\rvert^{2} + \lambda_{+}\chi_{\{u>0\}}+\lambda_{-}\chi_{\{u<0\}}\mathrm{d}m\\
			= & \int_{B_{1/10}^{k}(\eta_{k})} \lvert\nabla u_{k}\rvert^{2} + \frac{\rho_{k}^{2}}{\lvert u(y_{k})\rvert^{2}} \left( \lambda_{+}\chi_{\{u_{k}>0\}} + \lambda_{-}\chi_{\{u_{k}<0\}} \right) \mathrm{d}m_{k}.
		\end{split}
	\end{equation}
	In particular, $u_{k}$ minimizes the energy functional
	\begin{equation}
		E_{k}(v):=\int_{B_{1/10}^{k}(\eta_{k})} \lvert \nabla v\rvert^{2} + \lambda_{+}^{(k)}\chi_{\{v>0\}} + \lambda_{-}^{(k)}\chi_{\{v<0\}}\mathrm{d}m_{k},
	\end{equation}
	where
	\begin{equation}\label{equation_3_Lipschitz_continuity}
		\lambda_{\pm}^{(k)} := \frac{\rho_{k}^{2}}{\lvert u(y_{k})\rvert^{2}} \lambda_{\pm}\le \left(\frac{20}{k}\right)^{2}\lambda_{\pm}\to 0,\quad \text{as } k\to\infty.
	\end{equation}
	We may assume that $\lambda_{\pm}^{(k)}\le \lambda_{\pm}$. By Theorem \ref{theorem_1_Holder_continuity}, $u_{k}$ are equi-continuous in $B_{3/40}^{k}(\eta_{k})$. More precisely, there exists a constant $C>0$ such that
	\begin{equation}
		\lvert u_{k}(x)-u_{k}(y)\rvert \le Cd(x,y)^{\frac{2}{N+2}},\quad \forall x,y\in B_{3/40}^{k}(\eta_{k}),\forall k\in\mathbb{N}.
	\end{equation}
	Moreover, by \eqref{equation_1_Holder_continuity} in the proof of Theorem \ref{theorem_1_Holder_continuity}, it holds that
	\begin{equation}
		\int_{B_{3/40}^{k}(\eta_{k})} \lvert\nabla u_{k}\rvert^{2} \mathrm{d}m_{k} \le C_{K,N,\lambda_{\pm}} m_{k}\big(B_{3/40}^{k}(\eta_{k})\big).
	\end{equation}
	In particular, we have that
	\begin{equation}
		\limsup_{k\to\infty} \int_{B_{3/40}^{k}(\eta_{k})} \lvert \nabla u_{k}\rvert^{2}\mathrm{d}m_{k} \le C_{K,N,\lambda_{\pm}} \limsup_{k\to\infty} m_{k}\big(B_{3/40}^{k}(\eta_{k})\big)\le C_{K,N,\lambda_{\pm}} m_{\infty} \big(B_{1/10}^{\infty}(\eta_{\infty})\big)<\infty.
	\end{equation}
	
    For each $k\in\mathbb{N}$, let $\phi_{k}$ be a good cut-off function with $\phi_{k}=1$ on $B_{5/80}^{k}(\eta_{k})$ and $\mathop{\mathrm{supp}}\phi_{k} \subset B_{3/40}^{k}(\eta_{k})$. Recall that $\lvert \nabla \phi_{k}\rvert\le C_{K,N}$. In particular, $u_{k}\phi_{k}$ are uniformly bounded and equi-continuous. By the Arzel\'{a}-Ascoli theorem, $u_{k}\phi_{k}$ converges uniformly (hence $L^{2}$-strongly) to some $u_{\infty}:X_{\infty}\to\mathbb{R}$, up to a subsequence. Moreover, we have that $\lvert\nabla (u_{k}\phi_{k})\rvert^{2}\le 2u_{k}^{2}\lvert\nabla \phi_{k}\rvert^{2} + 2\phi_{k}^{2}\lvert\nabla u_{k}\rvert^{2}\le 2\big(\lvert\nabla \phi_{k}\rvert^{2} + \lvert\nabla u_{k}\rvert^{2}\big)$ and hence
    \begin{equation}
    	\sup_{k\in\mathbb{N}} \lvert \nabla (u_{k}\phi_{k})\rvert_{L^{2}} \le 2 \sup_{k\in\mathbb{N}} \lvert\nabla u_{k}\rvert_{L^{2}(B_{3/40}^{k}(\eta_{k}))} + 2\sup_{k\in\mathbb{N}} \lvert\nabla \phi_{k}\rvert_{L^{2}(B_{3/40}^{k}(\eta_{k}))}<\infty.
    \end{equation}
    Therefore, $u_{k}\phi_{k}$ $W^{1,2}$-weakly converges to $u_{\infty}$.
    
    \textbf{Step 4.} For each $f_{\infty}\in\mathop{\mathrm{Lip}_{0}}\big(B_{5/80}^{\infty}(\eta_{\infty})\big)$, there exists a sequence of $f_{k}\in\mathop{\mathrm{Lip}_{0}}\big(B_{5/80}^{k}(\eta_{k})\big)$ $W^{1,2}$-strongly converging to $f_{\infty}$ with $\sup\limits_{k} \lvert\nabla f_{k}\rvert_{L^{\infty}}<\infty$. Note that
    \begin{equation}
    	\sup_{k} \lvert \nabla(u_{k}\phi_{k})\cdot\nabla f_{k}\rvert_{L^{2}} \le \left( \sup_{k} \lvert\nabla f_{k}\rvert_{L^{\infty}}\right) \left(\sup_{k} \lvert \nabla(u_{k}\phi_{k})\rvert_{L^{2}}\right)<\infty.
    \end{equation}
    Therefore, we get that $\nabla (u_{k}\phi_{k})\cdot\nabla f_{k}$ $L^{2}$-weakly converges to $\nabla u_{\infty}\cdot\nabla f_{\infty}$. In particular, we get that
    \begin{equation}
    	\lim_{k\to\infty} \int_{B_{5/80}^{k}(\eta_{k})} \nabla (u_{k}\phi_{k})\cdot\nabla f_{k}\mathrm{d}m_{k}= \int_{B_{5/80}^{\infty}(\eta_{\infty})} \nabla u_{\infty} \cdot\nabla f_{\infty}\mathrm{d}m_{\infty}.
    \end{equation}
    By Lemma \ref{lemma_2_Lipschitz_continuity} and \eqref{equation_3_Lipschitz_continuity}, we get that
    \begin{equation}
    	- \int_{B_{5/80}^{k}(\eta_{k})} \nabla (u_{k}\phi_{k})\cdot\nabla f_{k}\mathrm{d}m_{k} \le \lvert \nabla f_{k}\rvert_{L^{2}}\sqrt{\lambda_{+}^{(k)}m_{k}\big(B_{5/80}^{k}(\eta_{k})\big)} \to 0,\quad\text{as } k\to\infty.
    \end{equation}
    Therefore, we get that
    \begin{equation}
    	- \int_{B_{5/80}^{\infty}(\eta_{\infty})} \nabla u_{\infty}\cdot\nabla f_{\infty}\mathrm{d}m_{\infty}\le 0.
    \end{equation}
    Since $f_{\infty}\in\mathop{\mathrm{Lip}_{0}}\big(B_{5/80}^{\infty}(\eta_{\infty})\big)$ is arbitrary, we get that $u_{\infty}$ is harmonic on $B_{5/80}^{\infty}(\eta_{\infty})$.
    
    Note that $u_{\infty}(\eta_{\infty})=0$ since $u_{k}(\eta_{k})=0$ and that $u_{\infty}\ge 0$ by \eqref{equation_2_Lipschitz_continuity}. By the strong maximum principle, we get that $u_{\infty}\equiv 0$ on $B_{5/80}^{\infty}(\eta_{\infty})$, contradicting to the last line in \eqref{equation_1_Lipschitz_continuity}.
\end{proof}
The main result of this section is the following theorem.
\begin{theorem}\label{theorem_1_Lipschitz_continuity}
	Let $u$ be a minimizer of $E$. If one of $u_{\pm}$ has linear growth locally on $\Omega$, then $u$ is locally Lipschitz continuous on $\Omega$. In particular, for the one-phase problem, $u$ is locally Lipschitz continuous on $\Omega$.
\end{theorem}
\begin{proof}
	Since $u$ is harmonic on the open sets $\{u_{\pm}>0\}$, $u$ is locally Lipschitz continuous on $\{u_{\pm}>0\}$. It remains to show that for each $B_{2r}(x_{0})\subset \Omega$ with $x_{0}\in\partial\{u>0\}\cup\partial\{u<0\}$, $u$ is Lipschitz continuous on $B_{r}(x_{0})$. Indeed, by Lemma \ref{lemma_3_Lipschitz_continuity}, there exists a constant $L>0$ such that
	\begin{equation}
		\lvert u\rvert(x)\le Ld(x,\{u=0\}),\quad\forall x\in B_{r}(x_{0}).
	\end{equation}
	Therefore, for each $x\in \{u_{\pm}>0\}\cap B_{r}(x_{0})$, we have that $u$ is harmonic in $B_{d(x,\{u=0\})}(x)$ and thus
	\begin{equation}
		\lvert \nabla u\rvert(x)\le C_{K,N,r}\frac{\lvert u\rvert(x)}{d(x,\{u=0\})} \le C_{K,N,r}L
	\end{equation}
	by Yau's gradient estimate. That is,
	\begin{equation}
		\lvert \nabla u\rvert_{L^{\infty}(B_{r}(x_{0}))} \le C_{K,N,r}L<\infty,
	\end{equation}
	which implies that $u$ is Lipschitz continuous on $B_{r}(x_{0})$.
\end{proof}
\begin{remark}\label{remark_1_Lipschitz_continuity}
	For the two-phase problem on $\mathbb{R}^{n}$, the original proof of local Lipchitz continuity of minimizers requires the Alt-Caffarelli-Friedman monotonicity formula \cite{Alt-Caffarelli-Friedman_1984}. On general $RCD(K,N)$-spaces, we lack a counterpart of the Alt-Caffarelli-Friedman monotonicity formula. This is the main difficulty in fully establishing the local Lipschitz continuity for the two-phase problem on $RCD(K,N)$-spaces.
\end{remark}

\section{Non-degeneracy}
In this section, we prove that minimizers are non-degenerate. We recall the Stokes formula that will be used later \cite[Remark 5.2]{Chan-Zhang-Zhu_2021}.
\begin{lemma}\label{lemma_1_Non-degeneracy}
	Let $(X,d,m)$ be an $RCD(K,N)$-space with $K\in\mathbb{R}$ and $N\in(1,\infty)$ and let $B_{R}(x_{0})\subset X$. Let $\phi\in C^{2}[0,R]$ and $\varphi:=\phi\circ d(\cdot,x_{0})$. If $\mathbf{\Delta}\varphi$ is a signed Radon measure, then for each $u\in W^{1,2}\big(B_{R}(x_{0})\big)\cap C\big(\overline{B}_{R}(x_{0})\big)$ it holds that 
	\begin{equation}
		\int_{B_{r}(x_{0})} u\mathbf{\Delta}\varphi = -\int_{B_{r}(x_{0})} \nabla u\cdot\nabla\varphi\mathrm{d}m + \phi^{\prime}(r)\frac{\mathrm{d}}{\mathrm{d}s}\bigg|_{s=r}\int_{B_{s}(x_{0})}u\mathrm{d}m
	\end{equation}
	for almost all $0<r<R$, and that
	\begin{equation}
		\begin{split}
			\int_{B_{r_{2}}(x_{0})\setminus B_{r_{1}}(x_{0})} u\mathbf{\Delta}\varphi = & - \int_{B_{r_{2}}(x_{0})\setminus B_{r_{1}}(x_{0})} \nabla u\cdot\nabla\varphi\mathrm{d}m\\
			& + \phi^{\prime}(r_{2})\frac{\mathrm{d}}{\mathrm{d}r}\bigg|_{s=r_{2}} \int_{B_{s}(x_{0})} u\mathrm{d}m - \phi^{\prime}(r_{1}) \frac{\mathrm{d}}{\mathrm{d}r}\bigg|_{s=r_{1}}\int_{B_{s}(x_{0})}u\mathrm{d}m
		\end{split}
	\end{equation}
	for almost all $0<r_{1}<r_{2}<R$.
\end{lemma}
The main result of this section is the following theorem.
\begin{theorem}\label{theorem_1_Non-degeneracy}
	Let $u$ be a minimizer of $E$. Then there exists a constant $C>0$, depending only on $K,N$, such that
	\begin{equation}
		\sup_{B_{R}(x_{0})} u_{\pm} \ge C\sqrt{\lambda_{\pm}}R
	\end{equation}
	for all $B_{R}(x_{0})\subset \Omega$ with $x_{0}\in \partial\{u_{\pm}>0\}$ and $R\in(0,1)$.
\end{theorem}
\begin{proof}
	We only prove the case of $x_{0}\in\partial\{u_{+}>0\}$. The case of $x_{0}\in\partial\{u_{-}>0\}$ is analogous. We define
	\begin{equation}
		\tilde{\phi}(t):=-\int_{t}^{1}\left( \frac{\sinh(s\sqrt{-\kappa})}{\sqrt{-\kappa}} \right)^{1-N} \mathrm{d}s,\quad s\in[0,1],
	\end{equation}
	where $\kappa:=\frac{\lvert K\rvert}{N-1}$. We define
	\begin{equation}
		\begin{split}
			& \phi(t):=\frac{\big(\tilde{\phi}(t)-\tilde{\phi}(1/2)\big)_{+}}{\tilde{\phi}(1)-\tilde{\phi}(1/2)},\\
			& \varphi(x):=\phi\left(\frac{d(x,x_{0})}{R}\right).
		\end{split}
	\end{equation}
	We define $M:=\sup\limits_{B_{R}(x_{0})}u_{+}$ and consider the function
	\begin{equation}
		v:=\begin{cases}
			u,\quad & \text{on } \Omega\setminus B_{R}(x_{0}),\\
			\min\{u,M\varphi\},\quad & \text{on } B_{R}(x_{0}).
		\end{cases}
	\end{equation}
	Note that $v\in\mathscr{A}_{g}$ and hence $E(u)\le E(v)$. A simple computation gives that
	\begin{equation}
		\int_{B_{R/2}(x_{0})} \lvert \nabla u_{+}\rvert^{2} + \lambda_{+}\chi_{\{u>0\}} \mathrm{d}m \le \int_{B_{R}(x_{0})\setminus B_{R/2}(x_{0})} \lvert \nabla v\rvert^{2} - \lvert \nabla u\rvert^{2} \mathrm{d}m.
	\end{equation}
	We define $w:=u-v$. Then $w=(u-M\varphi)_{+}$ on $B_{R}(x_{0})$ and hence
	\begin{equation}\label{equation_1_Non-degeneracy}
		\begin{split}
			\int_{B_{R/2}(x_{0})} \lvert \nabla u_{+}\rvert^{2} + \lambda_{+}\chi_{\{u>0\}}\mathrm{d}m \le & \int_{B_{R}(x_{0})\setminus B_{R/2}(x_{0})} -2\nabla w\cdot\nabla v - \lvert \nabla w\rvert^{2} \mathrm{d}m\\
			\le & -2M \int_{B_{R}(x_{0})\setminus B_{R/2}(x_{0})} \nabla w\cdot\nabla \varphi\mathrm{d}m.
		\end{split}
	\end{equation}
	By the Stokes formula, for almost all $\frac{R}{2}<r_{1}<r_{2}<R$ it holds that
	\begin{equation}
		\begin{split}
			& -\int_{B_{r_{2}}(x_{0})\setminus B_{r_{1}}(x_{0})} \nabla w\cdot\nabla\varphi\mathrm{d}m\\
			= & \int_{B_{r_{2}}(x_{0})\setminus B_{r_{1}}(x_{0})} w\mathbf{\Delta}\varphi + \frac{1}{R}\phi^{\prime}\left(\frac{r_{1}}{R}\right) \frac{\mathrm{d}}{\mathrm{d}s}\bigg|_{s=r_{1}} \int_{B_{s}(x_{0})} w\mathrm{d}m - \frac{1}{R}\phi^{\prime}\left(\frac{r_{2}}{R}\right) \frac{\mathrm{d}}{\mathrm{d}s}\bigg|_{s=r_{2}} \int_{B_{s}(x_{0})} w\mathrm{d}m.
		\end{split}
	\end{equation}
	Note that $w\ge 0$ and $\mathbf{\Delta}\varphi\le 0$ by the Laplacian comparison theorem for distance functions. Therefore
	\begin{equation}\label{equation_2_Non-degeneracy}
		-\int_{B_{r_{2}}(x_{0})\setminus B_{r_{1}}(x_{0})} \nabla w\cdot\nabla\varphi\mathrm{d}m \le  \frac{1}{R}\phi^{\prime}\left(\frac{r_{1}}{R}\right) \frac{\mathrm{d}}{\mathrm{d}s}\bigg|_{s=r_{1}} \int_{B_{s}(x_{0})} w\mathrm{d}m - \frac{1}{R}\phi^{\prime}\left(\frac{r_{2}}{R}\right) \frac{\mathrm{d}}{\mathrm{d}s}\bigg|_{s=r_{2}} \int_{B_{s}(x_{0})} w\mathrm{d}m.
	\end{equation}
	
	We estimate $\lim\limits_{r_{2}\to R}\frac{\mathrm{d}}{\mathrm{d}s}\bigg|_{s=r_{2}}\int_{B_{s}(x_{0})}w\mathrm{d}m$ as below. For each $\delta>0$, we have that
	\begin{equation}
		\frac{1}{\delta}\left\vert \int_{B_{r_{2}+\delta}(x_{0})} w\mathrm{d}m - \int_{B_{r_{2}}(x_{0})} w\mathrm{d}m \right\rvert \le \frac{m\left(B_{r_{2}+\delta}(x_{0})\setminus B_{r_{2}}(x_{0})\right)}{\delta} \sup_{B_{r_{2}+\delta}(x_{0}) \setminus B_{r_{2}}(x_{0})} \lvert w\rvert.
	\end{equation}
	By the Bishop-Gromov inequality \cite[Theorem 2.3]{Sturm_2006_2}, we have that (recall that $R/2<r_{2}<R$)
	\begin{equation}
		\limsup_{\delta \to 0} \frac{m\left(B_{r_{2}+\delta}(x_{0})\setminus B_{r_{2}}(x_{0})\right)}{\delta} \le C_{K,N,R} m\left(B_{r_{2}}(x_{0})\right) \le C_{K,N,R} m\left(B_{R}(x_{0})\right)<\infty.
	\end{equation}
	On the other hand, $w$ is continuous on the compact set $\overline{B}_{R}(x_{0})\setminus B_{R/2}(x_{0})$ and $w=0$ on $\partial B_{R}(x_{0})$. In particular, we have that
	\begin{equation}
		\lim_{\varepsilon\to 0} \sup_{B_{R}(x_{0})\setminus B_{R-\varepsilon}(x_{0})} \lvert w\rvert = 0.
	\end{equation}
	Hence, we have that
	\begin{equation}\label{equation_3_Non-degeneracy}
		\lim_{r_{2}\to R} \frac{\mathrm{d}}{\mathrm{d}s}\bigg|_{s=r_{2}}\int_{B_{s}(x_{0})}w\mathrm{d}m=0.
	\end{equation}
	Using the fact that $w=u_{+}$ on $\partial B_{R/2}(x_{0})$, we also have that
	\begin{equation}\label{equation_4_Non-degeneracy}
		\lim_{r_{1}\to R/2} \frac{\mathrm{d}}{\mathrm{d}s}\bigg|_{s=r_{1}}\int_{B_{s}(x_{0})} w\mathrm{d}m = \lim_{r_{1}\to R/2} \frac{\mathrm{d}}{\mathrm{d}s}\bigg|_{s=r_{1}}\int_{B_{s}(x_{0})} u_{+}\mathrm{d}m.
	\end{equation}
	
	Using the Stokes formula again, for almost all $r_{1}\in (R/2,R)$ we have that (recall that $0<R<1$)
	\begin{equation}\label{equation_5_Non-degeneracy}
		\begin{split}
			\frac{\mathrm{d}}{\mathrm{d}s}\bigg|_{s=r_{1}}\int_{B_{s}(x_{0})} u_{+}\mathrm{d}m = & \frac{1}{2r_{1}} \left( \int_{B_{r_{1}}(x_{0})}u_{+}\mathbf{\Delta}d(\cdot,x_{0})^{2} + \int_{B_{r_{1}}(x_{0})}\nabla u_{+}\cdot\nabla d(\cdot,x_{0})^{2}\mathrm{d}m \right)\\
			\le & \frac{C_{K,N}}{r_{1}} \int_{B_{r_{1}}(x_{0})} u_{+}\mathrm{d}m + \frac{1}{2r_{1}}\int_{B_{r_{1}}(x_{0})} 2d(\cdot,x_{0})\lvert\nabla u_{+}\rvert\mathrm{d}m\\
			\le & \frac{C_{K,N}}{r_{1}} \int_{B_{r_{1}}(x_{0})} M\chi_{\{u>0\}}\mathrm{d}m + \int_{B_{r_{1}}(x_{0})} \lvert \nabla u_{+}\rvert\mathrm{d}m\\
			\le & \frac{C_{K,N}}{r_{1}} \int_{B_{r_{1}}(x_{0})}M\chi_{\{u>0\}}\mathrm{d}m + \int_{B_{r_{1}}(x_{0})} \frac{1}{2} \left( \sqrt{\lambda_{+}}\chi_{\{u>0\}} + \frac{\lvert\nabla u_{+}\rvert^{2}}{\sqrt{\lambda_{+}}} \right) \mathrm{d}m\\
			\le & C_{K,N} \left( \frac{M}{\lambda_{+}r_{1}} + \frac{1}{\sqrt{\lambda_{+}}}\right) \int_{B_{r_{1}}(x_{0})} \lvert\nabla u_{+}\rvert^{2} + \lambda_{+}\chi_{\{u>0\}}\mathrm{d}m,
		\end{split}
	\end{equation}
	where we used the Laplacian comparison theorem for distance functions in the second line.
	
	By \eqref{equation_2_Non-degeneracy}, \eqref{equation_3_Non-degeneracy}, \eqref{equation_4_Non-degeneracy} and \eqref{equation_5_Non-degeneracy} and letting $r_{1}\to R/2$ and $r_{2}\to R$, we get that
	\begin{equation}\label{equation_6_Non-degeneracy}
		-\int_{B_{R}(x_{0})\setminus B_{R/2}(x_{0})} \nabla w\cdot\nabla \varphi\mathrm{d}m \le \frac{1}{R}\phi^{\prime}\left(\frac{1}{2}\right)C_{K,N}\left( \frac{M}{\lambda_{+}R} + \frac{1}{\sqrt{\lambda_{+}}}\right)\int_{B_{R/2}(x_{0})} \lvert\nabla u_{+}\rvert^{2} + \lambda_{+}\chi_{\{u>0\}}\mathrm{d}m.
	\end{equation}
	Since $x_{0}\in\partial\{u>0\}$, we have that $\int_{B_{R/2}(x_{0})} \lvert\nabla u_{+}\rvert^{2} + \lambda_{+}\chi_{\{u>0\}}\mathrm{d}m>0$. Hence, by \eqref{equation_1_Non-degeneracy} and \eqref{equation_6_Non-degeneracy}, we get that
	\begin{equation}
		\frac{2M}{R}\phi^{\prime}\left(\frac{1}{2}\right)C_{K,N}\left( \frac{M}{\lambda_{+}R}+ \frac{1}{\lambda_{+}}\right)\ge 1.
	\end{equation}
	Note that $\phi^{\prime}\left(\frac{1}{2}\right)$ depends only on $K,N$. Hence,
	\begin{equation}
		\frac{C_{K,N}M}{\sqrt{\lambda_{+}}R}\left( \frac{M}{\sqrt{\lambda_{+}}R}+ 1 \right)\ge 1.
	\end{equation}
	In particular, either $\frac{C_{K,N}M}{\sqrt{\lambda_{+}}R}\ge \frac{1}{2}$ or $\frac{M}{\sqrt{\lambda_{+}}R}+1\ge 2$, which implies that
	\begin{equation}
		\sup_{B_{R}(x_{0})} u_{+} = M \ge \min \left\{ \frac{\sqrt{\lambda_{+}}R}{2C_{K,N}}, \sqrt{\lambda_{+}}R \right\}.
	\end{equation}
\end{proof}

\section{Local finiteness of perimeter}
In this section, we prove that the free boundaries of minimizers has locally finite perimeter.
\begin{lemma}\label{lemma_1_Local_finiteness_of_perimeter}
	Let $u$ be a minimizer of $E$. Then for all $B_{2r}(x_{0})\subset \Omega$ with $x_{0}\in\partial\{u_{\pm}>0\}$ it holds that
	\begin{equation}
		\int_{B_{r}(x_{0})\cap\{0<u_{\pm}\le\varepsilon\}} \lvert\nabla u\rvert^{2}+\lambda_{\pm}\mathrm{d}m\le C\varepsilon,\quad\forall \varepsilon\in(0,\varepsilon),
	\end{equation}
	where the constant $C>0$ depends only on $K,N,r,\lvert\nabla u_{\pm}\rvert_{L^{2}(B_{2r}(x_{0}))}$ and $m\big(B_{2r}(x_{0})\big)$.
\end{lemma}
\begin{proof}
	We only prove the case of $x_{0}\in\partial\{u_{+}>0\}$. The case of $x_{0}\in\partial\{u_{-}>0\}$ is analogous.
	
	Let $\phi$ be a good cut-off function such that
	\begin{equation}
		\phi=1\quad\text{on } B_{r}(x_{0})\quad\text{and}\quad \phi=0\quad\text{on }X\setminus B_{2r}(x_{0}).
	\end{equation}
	Recall that $\lvert\nabla\phi\rvert\le C_{K,N,r}$.
	
	We define
	\begin{equation}
		v:=\begin{cases}
			u,\quad & \text{if } u\le 0,\\
			(1-\phi)u,\quad & \text{if } 0<u\le \varepsilon,\\
			\left(1-\frac{\varepsilon\phi}{u}\right)u,\quad & \text{if } u>\varepsilon.
		\end{cases}
	\end{equation}
	Then $v\in\mathscr{A}_{g}$ and hence $E(u)\le E(v)$. A simple computation gives that
	\begin{equation}\label{equation_1_Local_finiteness_of_perimeter}
		\begin{split}
			\int_{B_{2r}(x_{0})\cap\{u>0\}} \lvert\nabla u\rvert^{2} -\lvert\nabla v\rvert^{2}\mathrm{d}m \le & \int_{B_{2r}(x_{0})} \lambda_{+}\chi_{\{v>0\}} - \lambda_{+}\chi_{\{u>0\}}\mathrm{d}m\\
			\le & - \int_{B_{r}(x_{0})\cap\{0<u\le\varepsilon\}} \lambda_{+}\mathrm{d}m.
		\end{split} 
	\end{equation}
	
	We estimate $\int_{B_{2r}(x_{0})\cap\{u>0\}} \lvert\nabla u\rvert^{2} -\lvert\nabla v\rvert^{2}\mathrm{d}m$ as below.
	
	On $B_{2r}(x_{0})\cap\{0<u\le\varepsilon\}$, we have that $v=(1-\phi)u$ and therefore
	\begin{equation}\label{equation_2_Local_finiteness_of_perimeter}
		\begin{split}
			\int_{B_{2r}(x_{0})\cap\{0<u\le\varepsilon\}} \lvert\nabla u\rvert^{2} -\lvert\nabla v\rvert^{2}\mathrm{d}m = & \int_{B_{2r}(x_{0})\cap\{0<u\le\varepsilon\}} (2\phi-\phi^{2})\lvert\nabla u\rvert^{2} + 2(1-\phi)u\nabla u\cdot\nabla\phi - u^{2}\lvert\nabla\phi\rvert^{2}\mathrm{d}m\\
			\ge & \int_{B_{2r}(x_{0})\cap\{0<u\le\varepsilon\}} \chi_{B_{r}(x_{0})}\lvert\nabla u\rvert^{2} - 2\varepsilon\lvert\nabla u\rvert\lvert\nabla\phi\rvert - \varepsilon^{2}\lvert\nabla \phi\rvert^{2}\mathrm{d}m\\
			\ge & \int_{B_{r}(x_{0})\cap\{0<u\le\varepsilon\}} \lvert\nabla u\rvert^{2}\mathrm{d}m - 2\varepsilon \lvert \nabla u_{+}\rvert_{L^{2}(B_{2r}(x_{0}))} \lvert\nabla\phi\rvert_{L^{2}} - \varepsilon^{2}\lvert\nabla\phi\rvert^{2}_{L^{2}}\\
			\ge & \int_{B_{r}\cap\{0<u\le\varepsilon\}} \lvert\nabla u\rvert^{2}\mathrm{d}m -C_{1}\varepsilon,
		\end{split}
	\end{equation}
	where the constant $C_{1}>0$ depends only on $K,N,r,\lvert\nabla u_{+}\rvert_{L^{2}(B_{2r}(x_{0}))}$ and $m\big(B_{2r}(x_{0})\big)$.
	
	On $B_{2r}(x_{0})\cap\{u>\varepsilon\}$, we have that $v=u-\varepsilon\phi$ and therefore
	\begin{equation}\label{equation_3_Local_finiteness_of_perimeter}
		\begin{split}
			\int_{B_{2r}(x_{0})\cap\{u>\varepsilon\}} \lvert\nabla u\rvert^{2} -\lvert\nabla v\rvert^{2}\mathrm{d}m = & \int_{B_{2r}(x_{0})\cap\{u>\varepsilon\}} 2\varepsilon\nabla u\cdot\nabla\phi - \varepsilon^{2}\lvert\nabla\phi\rvert^{2}\mathrm{d}m\\
			\ge & -2\varepsilon \lvert u_{+}\rvert_{L^{2}(B_{2r}(x_{0}))} \lvert\nabla\phi\rvert_{L^{2}} - \varepsilon^{2}\lvert\nabla\phi\rvert_{L^{2}}^{2}\\
			\ge & -C_{2}\varepsilon
		\end{split}
	\end{equation}
	where the constant $C_{2}>0$ depends only on $K,N,r,\lvert\nabla u_{+}\rvert_{L^{2}(B_{2r}(x_{0}))}$ and $m\big(B_{2r}(x_{0})\big)$.
	
	By \eqref{equation_1_Local_finiteness_of_perimeter}, \eqref{equation_2_Local_finiteness_of_perimeter} and \eqref{equation_3_Local_finiteness_of_perimeter}, we get that
	\begin{equation}
		\int_{B_{r}(x_{0})\cap\{0<u\le\varepsilon\}} \lvert\nabla u\rvert^{2}+\lambda_{+}\mathrm{d}m \le (C_{1}+C_{2})\varepsilon.
	\end{equation}
\end{proof}

The main result of this section is the following theorem.
\begin{theorem}\label{theorem_1_Local_finiteness_of_perimeter}
	Let $u$ be a minimizer of $E$. Then the free boundaries $\partial\{u_{\pm}>0\}\cap\Omega$ has locally finite perimeter in $\Omega$. More precisely, for each $B_{2r}(x_{0})\subset\Omega$ with $x_{0}\in\partial\{u_{\pm}>0\}$ it holds that
	\begin{equation}
		\mathop{\mathrm{Per}}\left(\{u_{\pm}>0\};B_{r}(x_{0})\right) \le C,
	\end{equation}
	where the constant $C>0$ depends only on $K,N,r,\lambda_{\pm},\lvert\nabla u_{\pm}\rvert_{L^{2}(B_{2r}(x_{0}))}$ and $m\big(B_{2r}(x_{0})\big)$.
\end{theorem}
\begin{proof}
	Let $\varepsilon\in(0,1)$. By the coarea formula, it holds that
	\begin{equation}
		\int_{\varepsilon/2}^{\varepsilon} \mathop{\mathrm{Per}}\left(\{u_{\pm}>t\};B_{r}(x_{0})\right)\mathrm{d}t = \int_{\{\varepsilon/2\le u_{\pm}<\varepsilon\}} \chi_{B_{r}(x_{0})} \mathrm{d}\lvert Du\rvert.
	\end{equation}
	Note that
	\begin{equation}
		\begin{split}
			\int_{\{\varepsilon/2\le u_{\pm}<\varepsilon\}} \chi_{B_{r}(x_{0})} \mathrm{d}\lvert Du\rvert = & \int_{B_{r}(x_{0})\cap\{\varepsilon/2\le u_{\pm}< \varepsilon\}} \lvert\nabla u\rvert\mathrm{d}m\\
			\le & \int_{B_{r}(x_{0})\cap\{0<u_{\pm}\le\varepsilon\}} \frac{\lvert\nabla u\rvert^{2} + 1}{2}\mathrm{d}m.
		\end{split}
	\end{equation}
	By Lemma \ref{lemma_1_Local_finiteness_of_perimeter}, there exists a constant $C>0$, depending only on $K,N,r,\lambda_{\pm},\lvert\nabla u_{\pm}\rvert_{L^{2}(B_{2r}(x_{0}))}$ and $m\big(B_{2r}(x_{0})\big)$, such that
	\begin{equation}
		\int_{B_{r}(x_{0})\cap\{0<u_{\pm}\le\varepsilon\}} \frac{\lvert\nabla u\rvert^{2} + 1}{2}\mathrm{d}m \le C\varepsilon.
	\end{equation}
	Therefore, we get that
	\begin{equation}
		\int_{\varepsilon/2}^{\varepsilon} \mathop{\mathrm{Per}}\left(\{u_{\pm}>t\};B_{r}(x_{0})\right)\mathrm{d}t \le C\varepsilon.
	\end{equation}
	In particular, there exists a $t_{\varepsilon}\in (\varepsilon/2,\varepsilon)$ such that
	\begin{equation}
		\mathop{\mathrm{Per}}\left(\{u_{\pm}>t_{\varepsilon}\};B_{r}(x_{0})\right) \le \frac{2}{\varepsilon}C\varepsilon = 2C.
	\end{equation}
	By the lower semi-continuity of perimeter measure, we get that
	\begin{equation}
		\mathop{\mathrm{Per}}\left(\{u_{\pm}>0\};B_{r}(x_{0})\right) \le \liminf_{\varepsilon\to 0} \mathop{\mathrm{Per}} \left(\{u_{\pm}>t_{\varepsilon};B_{r}(x_{0})\}\right) \le 2C.
	\end{equation}
\end{proof}
\begin{remark}\label{remark_1_Local_finiteness_of_perimeter}
	In the proof of Theorem \ref{theorem_1_Local_finiteness_of_perimeter}, we use the fact that the function $v:=\min\{\max\{u_{\pm},\varepsilon/2\},\varepsilon\}$ is locally Lipschitz continuous on $\Omega$. Hence $\lvert Dv\rvert = \lvert \nabla v\rvert m$. Note that $\lvert \nabla v\rvert=\lvert\nabla u_{\pm}\rvert$ on $\{\varepsilon/2\le u_{\pm}\le\varepsilon\}$.
\end{remark}

\section*{Acknowledgments}
	
	\bibliographystyle{alpha}
	\bibliography{biblio}
	
\end{document}